\documentclass{article}
\usepackage[utf8]{inputenc}

\title{Porous numbers}
\author{R\"udiger Jehn }
\date{April 2021}

\usepackage[colorlinks,citecolor=blue,urlcolor=blue,bookmarks=false,hypertexnames=true]{hyperref} 
\begin{document}

\maketitle

{\bf Abstract}\\

The concept of porous numbers is presented. A number $k$ which is not a multiple of 10 is called {\it porous} if every number $m$ with sum of digits = $k$ and $k$ a divisor of both $m$ and digit reversal of $m$ has a zero in its digits. It is proved that 11, 37, 74, 101 and 121 are the only porous numbers smaller than 1000.

\section{Introduction}
A number $k$ which is not a multiple of 10 is called {\it porous} if every number $m$ that fulfills these 3 requirements:
\begin{enumerate}
    \item $k \; \vert \; m$
    \item $k \; \vert \; rev(m)$
    \item the sum of the digits of $m = k$ 
\end{enumerate}
must have at least one digit which is a zero. rev($m$) is the number m reversed (e.g.\ m = 123, then rev($m$) = 321). 
The numbers are called porous because when you "open" them with any $m$, there are always voids in $m$.\\

\section{List of non-porous numbers}

By definition if $k$ is a multiple of 10 it is non-porous. In fact, for any $m$, rev($m$) has a number different from 0 at the end and hence is not divisible by $k$ and no $m$ exists at all that could fulfill the 3 requirements.\\

For all other numbers a computer program was executed to find an $m$ not containing any zeros and fulfilling the 3 requirements. For $k$ smaller than 200, the results of the sequence A333666 \cite{A333666} were taken if $m$ did not contain any zero. If they had a zero, the search was extended until an $m$ without a zero was found. For $k$ larger than 200 most of the $m$ were found by concatenating palindromes of length up to 10. For example: $p_1 = 6498946$ is a palindrome which fulfills $202 \; \vert \; p$ and since it is a palindrome it also fulfills $202 \; \vert \; rev(p)$. The sum of the digits of $p_1$ is 46. $p_2 = 25452$ is a second palindrome that also fulfills the same requirements. The sum of the digits of $p_2$ is 18. Since 202 = 4 x 46 + 18, we concatenate 4 times $p_1$ with $p_2$ and the result 649894664989466498946649894625452 is the $m$ we were looking for.\\

The search was automatised and for most numbers an $m$ was quickly found except for multiples of 37, 101 and 121 which still had to be treated "manually". A file containing an $m$ for all non-porous numbers for $k$ up to 1000 was produced \cite{A337832}.

\section{List of porous numbers}

A brute force method to prove that a number is porous is to test all possible numbers $m$ that have a sum of digits $k$ and that do not contain a zero if both $k$ and rev($k$) are divisible by $m$. For $k = 11$ there are just 1021 numbers to be tested (1021 = A104144(19), there is an offset of 8 in this sequence \cite{A104144}). For $k = 37$ we get already 66 billion candidates, which is still doable, but at the latest for $k = 74$ the number of possibilities of $8.85e21$ is exceeding normal computing powers. Instead analytical reasoning will bring the proofs.

\subsection{Proof that 11 is a porous number}

Let "$m_{s-1} ... m_3 m_2 m_1 m_0$" be a number $m$ with $s$ digits that fulfills for $k=11$ the three requirements listed in the introduction. We define:\\ 

A = $m_0 + m_2 + m_4 + $ ... ,

B = $m_1 + m_3 + m_5 + $ ...  \\

A divisibility rule for 11 requires that the alternating sum of the digits must be divisible by 11. Hence:

\begin{equation}
   A - B = j \cdot 11
      \label{eq:1}
\end{equation}

Since the sum of the digits is 11, we have

\begin{equation}
   A + B = 11
    \label{eq:2}
\end{equation}

Adding eq.~\ref{eq:1} and \ref{eq:2} yields
$$ 2 A = (j + 1) \cdot 11 $$ therefore A must be 0 or 11. If A is 11, then B is 0. This means either A or B must be zero and $m$ must contain a zero. Hence 11 is a porous number.

\subsection{Proof that 37 is a porous number}

Let "$m_{s-1} ... m_3 m_2 m_1 m_0$" be a number $m$ with $s$ digits that fulfills for $k=37$ the three requirements listed in the introduction. $m_{s-1}$ is not allowed to be zero. We define:\\ 

A = $m_0 + m_3 + m_6 + $ ... ,

B = $m_1 + m_4 + m_7 + $ ... and 

C = $m_2 + m_5 + m_8 + $ ...\\

When $10^i$ ($i \ge 0)$ is divided by 37, the remainder is either 1, 10 or 26. We define:
$$
\delta_i = \left\{
    \begin{array}{rl}
        1 & \mbox{if } mod(i,3) = 0  \\
       10 & \mbox{if } mod(i,3) = 1  \\
       26 & \mbox{if } mod(i,3) = 2  \\
    \end{array}
\right.
$$

With this definition the powers of 10 can be written as:\\

$10^i = \alpha_i \cdot 37 + \delta_i$\\

And $m$ can be written as:

$$m = \sum_{i=0}^{s-1} m_i 10^i =  \sum_{i=0}^{s-1} m_i (\alpha_i \cdot 37 + \delta_i) $$\\

Since $37 \; \vert \; m$ it follows
\begin{equation}
   \sum_{i=0}^{s-1} m_i \delta_i  = l_1 \cdot 37
\end{equation}

Hence we get:

\begin{equation}
   A + 10 B + 26 C = l_1 \cdot 37 
   \label{eq:37}
\end{equation}

Eq.~\ref{eq:37} written for the reversed number of $m$ reads:

$$
\begin{array}{rl}
    A + 10 C + 26 B = j \cdot 37 & \mbox{if } mod(s,3) = 0 \mbox{, e.g. m has 13 digits} \\
    B + 10 A + 26 C = j \cdot 37 & \mbox{if } mod(s,3) = 1 \\
    26 A + 10 B + C = j \cdot 37 & \mbox{if } mod(s,3) = 2 \\
    \end{array} 
$$

For all three cases of $s$ it is easy to demonstrate that A, B and C must be multiples of 37. We show it here only for the second case with $mod(s,3) = 1$. Subtracting the equation for the reversed number from the equation for the original number gives
$$9 (B - A) = (l_1 - j) \cdot 37 $$
Therefore $B = A + l_2 \cdot 37$. Since $C = 37 - A - B$ we get $C = - 2 A + (1 - l_2) \cdot 37$. Inserting the expressions for $B$ and $C$ into eq.~\ref{eq:37} gives:

$$ 41 A = 37 (26 - l_1 - 16 l_2) $$
which is only possible if $A$ is a multiple of 37. Then $B$ must also be a multiple of 37 and finally also $C$. The only possibility to distribute 37 over three numbers which are all multiples of 37 is to allocate the whole 37 to one of the three and the other two must be zero. Therefore all $m$ must have at least 8 digits with a zero (see e.g.~A333666(37) = 1,009,009,009,009 \cite{A333666}).  \\

\subsection{Proof that 74 is a porous number}

This proof is very similar to the proof for $k = 37$. This time \\

$10^i = \alpha_i \cdot 74 + \gamma_i + \delta_i$\\

with 
$$
\gamma_i = \left\{
    \begin{array}{rl}
        -37 & \mbox{if } i = 0  \\
          0 & \mbox{else } \\
    \end{array}
\right.
$$ 
and 
$$
\delta_i = \left\{
    \begin{array}{rl}
        38 & \mbox{if } mod(i,3) = 0  \\
        10 & \mbox{if } mod(i,3) = 1  \\
        26 & \mbox{if } mod(i,3) = 2  \\
    \end{array}
\right.
$$

Now $m$ can be written as:

$$m = \sum_{i=0}^{s-1} m_i 10^i =  \sum_{i=0}^{s-1} m_i (\alpha_i \cdot 74 + \gamma_i  + \delta_i) $$\\

Since $74 \; \vert \; m$ it follows
\begin{equation}
   \sum_{i=0}^{s-1} m_i (\gamma_i  + \delta_i)  = l_1 \cdot 74
\end{equation}

Multiplying the individual terms gives
$$ 38 A + 10 B + 26 C - 37 m_0 = l_1 \cdot 74$$

$m_0$ must be an even number because $74 \; \vert \; m$ and therefore $37 m_0$ is a multiple of 74.\\

Hence we get:
$$ 38 A + 10 B + 26 C = l_2 \cdot 74$$

or 
\begin{equation}
   19 A + 5 B + 13 C = l_2 \cdot 37
   \label{eq:74}
\end{equation}

Eq.~\ref{eq:74} written for the reversed number of $m$ reads:

$$
\begin{array}{rl}
    19 A + 5 C + 13 B = j \cdot 37 & \mbox{if } mod(s,3) = 0 \mbox{, e.g. m has 16 digits} \\
    19 B + 5 A + 13 C = j \cdot 37 & \mbox{if } mod(s,3) = 1 \\
    13 A + 5 B + 19 C = j \cdot 37 & \mbox{if } mod(s,3) = 2 \\
    \end{array} 
$$

For all three cases of $s$ it is straightforward to demonstrate that $A$, $B$ and $C$ must be multiples of 37. Since the sum of the three is 74, at least one of them must be zero and therefore $m$ must have digits with zeros.

\subsection{Proof that 101 is a porous number}

For this proof we will exploit the fact that $101 \; \vert \; m$ if the "alternating sum of blocks of two" divides $m$.
As above, $m =$ "$m_{s-1} ... m_3 m_2 m_1 m_0$" is a number consisting of $s$ digits that fulfills for $k=101$ the three requirements listed in the introduction. We define:\\ 

A = $m_0 - m_2 + m_4 - m_6 + $ ... ,

B = $m_1 - m_3 + m_5 - m_7 + $ ...  \\

With this definition, the alternating sum of blocks of two, which we will call {\it alterdigitsum2,} is $A + 10 B$\\

From the divisibility rule for 101 it follows: 
\begin{equation}
   A + 10 B = l_1 \cdot 101
   \label{eq:101_1}
\end{equation}

Depending on the number of digits $s$, the alterdigitsum2 of $rev(m)$ will be

$$
\begin{array}{rl}
     - 10 A - B & \mbox{if } mod(s,4) = 0 \mbox{, e.g. m has 12 digits} \\
     A - 10 B & \mbox{if } mod(s,4) = 1 \\
     10 A + B & \mbox{if } mod(s,4) = 2 \\
     -A + 10 B & \mbox{if } mod(s,4) = 3 \\
\end{array} 
$$

If $s$ is odd then alterdigitsum2(rev(m)) = $ \pm (A - 10 B)$\\

From the divisibility rule for 101 it follows: 
\begin{equation}
   A - 10 B = l_2 \cdot 101
   \label{eq:101_2}
\end{equation}

Adding eqs.~\ref{eq:101_1} and \ref{eq:101_2} yields:

$$ 2 A = (l_1 + l_2) \cdot 101 $$
which is only possible if A is a multiple of 101. And if $A$ is a multiple of 101 also $B$ must be a multiple of 101.\\

If $s$ is even then alterdigitsum2(rev(m)) = $ \pm (10 A + B)$\\

From the divisibility rule for 101 it follows: 
\begin{equation}
   10 A + B = l_3 \cdot 101
   \label{eq:101_3}
\end{equation}

Subtracting eq.~\ref{eq:101_1} from 10 times eq.~\ref{eq:101_3} yields:

$$ 99 A = (10 l_3 - l_1) \cdot 101 $$
which is only possible if A is a multiple of 101. And if $A$ is a multiple of 101 also $B$ must be a multiple of 101.\\

What are the options then? $A$ and $B$ cannot be both zero because the sum of all digits, 101, is an odd number which means that if $A$ is even then $B$ must be odd and vice versa. 
Hence either $A$ or $B$ is $\pm 101$. But then all non-zero digits must be in the positions $j, j+4, j+8, ...$ with $j$ a number from 0 to 3 such that the $m_j$ can add up to 101. All other digits must be zero and therefore 101 is a porous number.\\

\subsection{Proof that 121 is a porous number}

As above, $m =$ "$m_{s-1} ... m_3 m_2 m_1 m_0$" is a number consisting of $s$ digits that fulfills for $k=121$ the three requirements listed in the introduction.
$s \ge 14$ because even if we fill 13 digits with a "9" the sum of digits is only 117. We define:\\ 

A = $m_0 + m_2 + m_4 + $ ... ,

B = $m_1 + m_3 + m_5 + $ ...  \\

Since $121 \; \vert \; m$ also $11 \; \vert \; m$ and a divisibility rule for number 11 requires that the alternating sum of the digits must be divisible by 11. Hence:

\begin{equation}
   A - B = j \cdot 11
    \label{eq:121_1}
\end{equation}

Since the sum of all digits is 121 we also know:

\begin{equation}
    A + B = 121 = 11 \cdot 11
    \label{eq:121_2}
\end{equation}

Adding eqs.~\ref{eq:121_1} and \ref{eq:121_2} gives:

\begin{equation}
2 A  = (j + 11) \cdot 11
\end{equation}

which means $A$ must be a multiple of 11, i.e. $j_1 \cdot 11$ and therefore $B = (11 - j_1) \cdot 11$. \\

Both $A$ and $B$ must be multiples of 11 and $j$ must be an odd number. \\

Next we exploit the fact that the multiplicative order of 10 modulo 121 equals 22:\\

$10^i = \alpha_i \cdot 121 + \beta_i \cdot 11 - 1^{i+1}$\\

with $\beta = (0, 1, 9, 3, 7, 5, 5, 7, 3, 9, 1, 0, 10, 2, 8, 4, 6, 6, 4, 8, 2, 10)$. After 22 numbers the sequence repeats, i.e. $\beta_i$ = $\beta_{i + 22}$.\\

Now $m$ can be written as:

$$m = \sum_{i=0}^{s-1} m_i 10^i =  \sum_{i=0}^{s-1} m_i (\alpha_i \cdot 121 + \beta_i \cdot 11 - 1^{i+1}) $$\\

Since $121 \; \vert \; m$ it follows
\begin{equation}
   \sum_{i=0}^{s-1} m_i ( \beta_i \cdot 11 - 1 ^ {i+1})  = l_1 \cdot 121
    \label{eq:121_4}
\end{equation}

Note, $\sum_{i=0}^{s-1}  - 1 ^ {i+1} m_i$ is the alternating sum of the digits, i.~e.~$A - B$. The first part of the sum is rewritten as:\\

$11 \cdot \left [  \vec{A} \cdot \beta_A + \vec{B} \cdot \beta_B \right ]$ \\

where $\vec{A} \cdot \beta_A = m_0 \beta_0 + m_2 \beta_2 + m_4 \beta_4 + ...$ and $\vec{B} \cdot \beta_B = m_1 \beta_1 + m_3 \beta_3 + m_5 \beta_5 + ...$\\

With this notation we can rewrite eq.~\ref{eq:121_4} as
\begin{equation}
   11 \cdot \left [ \vec{A} \cdot \beta_A + \vec{B} \cdot \beta_B \right ] + A - B = l_1 \cdot 121
   \label{eq:121_5}
\end{equation}

\subsubsection{Number of digits of $m$ is odd}

Now we formulate eq.~\ref{eq:121_5} for the reversed number $m$. We only need to replace $\beta_0$ with $\beta_{s-1}$, $\beta_1$ with $\beta_{s-2}$ and so forth. The result is:
\begin{equation}
   11 \cdot \left [ \vec{A} \cdot \bar{\beta_A} + \vec{B} \cdot  \bar{\beta_B} \right ] + A - B = l_2 \cdot 121
   \label{eq:121_rev_odd}
\end{equation}
with $\bar{\beta_A} = (\beta_{s-1}, \beta_{s-3}, \beta_{s-5}, ...)$\\

Adding eqs.~\ref{eq:121_5} and \ref{eq:121_rev_odd} yields:
\begin{equation}
   11 \cdot \left [ \vec{A} \cdot (\beta_A + \bar{\beta_A}) + \vec{B} \cdot  (\beta_B + \bar{\beta_B}) \right ] + 2(A - B) = (l_1 + l_2) \cdot 121
\label{eq:121_betas}
\end{equation}

Now we have to evaluate eq.~\ref{eq:121_betas} for s = 15, 17, ... until 35. Afterwards, starting from $s = 37$, the results are repeating because $\beta_i$ = $\beta_{i + 22}$. We notice that all elements of $\beta_A + \bar{\beta_A}$ have the same remainder $r$ when divided by 11 and $\beta_B + \bar{\beta_B}$ all have the same remainder $11 - r$ when divided by 11. In Table~\ref{table:1} the remainders~$r$ are printed for the number of digits $s$ before they repeat: \\

\begin{table}[h!]
\centering
\begin{tabular}[h]{|c|c|c|c|c|c|c|c|c|c|c|c|}
\hline
s & 15 & 17 & 19 & 21 & 23 & 25 & 27 & 29 & 31 & 33 & 35 \\
\hline
$\beta_A + \bar{\beta_A}$ & 8 & 6 & 4 & 2 & 0 & 9 & 7 & 5 & 3 & 1 & 10 \\
$\beta_B + \bar{\beta_B}$ & 3 & 5 & 7 & 9 & 0 & 2 & 4 & 6 & 8 & 10 & 1 \\
\hline 
\end{tabular}
\caption{Remainders of all elements of $\beta_A + \bar{\beta_A}$ and $\beta_B + \bar{\beta_B}$ divided by 11}
\label{table:1}
\end{table}

With this information eq.~\ref{eq:121_betas} can be written as:

$$ 11 \cdot \left [ r \cdot A + (11 - r ) \cdot B + l_3 \cdot 11 \right] + 2(A - B) = (l_1 + l_2) \cdot 121 $$

Since $A$ and $B$ are both multiples of 11 this means that $11 \cdot \left [ r \cdot A + (11 - r ) \cdot B + l_3 \cdot 11 \right ]$ is a multiple of 121 and therefore

$$ 2(A - B) = l_4 \cdot 121 $$

Since $A - B = j \cdot 11$ this equation can only be fulfilled if $j$ itself is a multiple of 11.

\subsubsection{Number of digits of $m$ is even}

If $s$ is even, the alternating digit sum of rev($m$) has the opposite sign and therefore we get instead of eq.~\ref{eq:121_rev_odd} this equation:
\begin{equation}
   11 \cdot \left [ \vec{A} \cdot \bar{\beta_A} + \vec{B} \cdot  \bar{\beta_B} \right ] + B - A = l_5 \cdot 121
    \label{eq:121_rev_even}
\end{equation}

In this case we subtract eq.~\ref{eq:121_rev_even} from eq.~\ref{eq:121_5}. This yields:
\begin{equation}
   11 \cdot \left [ \vec{A} \cdot (\beta_A - \bar{\beta_A}) + \vec{B} \cdot  (\beta_B - \bar{\beta_B}) \right ] + 2(A - B) = (l_1 - l_5) \cdot 121
\label{eq:121_betas_even}
\end{equation}

Now we have to evaluate eq.~\ref{eq:121_betas_even} for s = 14, 16, ... 34. Also in this case all elements of $\beta_A - \bar{\beta_A}$ have the same remainder $r$ when divided by 11 and $\beta_B - \bar{\beta_B}$ all have the same remainder $11 - r$ when divided by 11. In Table~\ref{table:2} the remainders~$r$ are printed for the number of digits $s$ before they repeat: \\

\begin{table}[h!]
\centering
\begin{tabular}[h]{|c|c|c|c|c|c|c|c|c|c|c|c|}
\hline
s & 14 & 16 & 18 & 20 & 22 & 24 & 26 & 28 & 30 & 32 & 34 \\
\hline
$\beta_A - \bar{\beta_A}$ & 9 & 7 & 5 & 3 & 1 & 10 & 8 & 6 & 4 & 2 & 0\\
$\beta_B - \bar{\beta_B}$ & 2 & 4 & 6 & 8 & 10 & 1 & 3 & 5 & 7 & 9 & 0 \\
\hline 
\end{tabular}
\caption{Remainders of all elements of $\beta_A - \bar{\beta_A}$ and $\beta_B - \bar{\beta_B}$ divided by 11}
\label{table:2}
\end{table}

With this information eq.~\ref{eq:121_betas_even} can be written as:

$$ 11 \cdot \left [ r \cdot A + (11 - r ) \cdot B + l_6 \cdot 11 \right] + 2(A - B) = (l_1 - l_5) \cdot 121 $$

Since $A$ and $B$ are both multiples of 11 this means that $11 \cdot \left [ r \cdot A + (11 - r ) \cdot B + l_6 \cdot 11 \right ]$ is a multiple of 121 and therefore

$$ 2(A - B) = l_7 \cdot 121 $$

Since $A - B = j \cdot 11$ this equation can only be fulfilled if $j$ itself is a multiple of 11. 

\subsubsection{Final conclusion for $k = 121$}

For both $s$ even and odd, it was shown that $j$ must be a multiple of 11.
Since $j$ is not zero we need to have $A - B = \pm 121$, i.e. either $A$ or $B$ must be 0 and the other must be 121. Hence 121 is a porous number. \\

\section{Outlook}

The steadily increasing number of possibilities to construct an $m$ for a given $k$ suggests that the 5 numbers which were proven in the previous chapter might be the only porous numbers. But a mathematical proof for this conjecture seems a big challenge.

\bibliography{references} 
\bibliographystyle{plain} 

\end{document}